\documentclass[12pt,reqno]{amsart}
\usepackage{amsmath,amssymb,amsfonts,amsthm,url,xspace,graphicx,psfrag}
\usepackage[pdftitle={The approach to criticality in sandpiles},
            pdfauthor={Anne Fey, Lionel Levine, David B. Wilson},
            bookmarks=true,bookmarksopen=true,bookmarksopenlevel=2,
            hyperfootnotes=false]{hyperref}
\usepackage[all]{hypcap}
\usepackage{multirow}
\usepackage[active]{srcltx}

\newcommand{\sref}[1]{\hyperref[#1]{\S~\ref*{#1}}}
\newcommand{\aref}[1]{\hyperref[#1]{Appendix~\ref*{#1}}}
\newcommand{\lref}[1]{\hyperref[#1]{Lemma~\ref*{#1}}}
\newcommand{\tref}[1]{\hyperref[#1]{Theorem~\ref*{#1}}}
\newcommand{\cref}[1]{\hyperref[#1]{Corollary~\ref*{#1}}}
\newcommand{\fref}[1]{\hyperref[#1]{Figure~\ref*{#1}}}
\newcommand{\pref}[1]{\hyperref[#1]{Proposition~\ref*{#1}}}
\renewcommand{\sref}[1]{\hyperref[#1]{section~\ref*{#1}}}

\newcommand{\old}[1]{}
\newcommand{\ceq}{\stackrel{?}{=}}

\def\clap#1{\hbox to 0pt{\hss#1\hss}}

\def\int{\mbox{int}}

 \renewcommand{\MRhref}[2]{\href{http://www.ams.org/mathscinet-getitem?mr=#1}{MR#2}}

\makeatletter
\def\@strippedMR{}
\def\@scanforMR#1#2#3\endscan{%
  \ifx#1M\ifx#2R\def\@strippedMR{#3}%
  \else\def\@strippedMR{#1#2#3}%
  \fi\fi}
\def\@rst #1 #2other{#1}
\renewcommand\MR[1]{\relax\ifhmode\unskip\spacefactor3000 \space\fi
  \@scanforMR#1\endscan
  \MRhref{\expandafter\@rst \@strippedMR other}{\@strippedMR}}
\newcommand\MRs[1]{\relax\ifhmode\unskip\spacefactor3000 \space\fi
  \@scanforMR#1\endscan
  \MRhref{\@strippedMR}{\@strippedMR}}
\makeatother

\newcommand{\floor}[1]{\left\lfloor {#1} \right\rfloor}
\newcommand{\ceiling}[1]{\left\lceil {#1} \right\rceil}

\newtheorem{theorem}{Theorem}[section]
\newtheorem{prop}[theorem]{Proposition}
\newtheorem{lemma}[theorem]{Lemma}
\newtheorem{corollary}[theorem]{Corollary}
\newtheorem{conjecture}[theorem]{Conjecture}

\theoremstyle{remark}

\theoremstyle{definition}

\def\mod{\mbox{mod }}

\DeclareSymbolFont{AMSb}{U}{msb}{m}{n}
\DeclareMathSymbol{\C}{\mathbin}{AMSb}{"43}
\DeclareMathSymbol{\E}{\mathbin}{AMSb}{"45}
\DeclareMathSymbol{\EE}{\mathbin}{AMSb}{"45}
\DeclareMathSymbol{\N}{\mathbin}{AMSb}{"4E}
\DeclareMathSymbol{\PP}{\mathbin}{AMSb}{"50}
\DeclareMathSymbol{\Q}{\mathbin}{AMSb}{"51}
\DeclareMathSymbol{\R}{\mathbin}{AMSb}{"52}
\DeclareMathSymbol{\Z}{\mathbin}{AMSb}{"5A}
\newcommand{\one}{{\mathbf 1}}

\newcommand{\odd}{\text{odd}}

\def\rcs $#1: #2 ${\expandafter\def\csname rcs#1\endcsname {#2}}
\rcs $Date: 2009/11/05 03:47:00 $

\begin{document}
\title{The approach to criticality in sandpiles}
\author{Anne Fey \and Lionel Levine \and David B. Wilson}

\address{Anne Fey, Delft Institute of Applied Mathematics, Delft University of Technology, The Netherlands, \texttt{\url{http://dutiosc.twi.tudelft.nl/~anne}}}
\address{Lionel Levine, Department of Mathematics, Massachusetts Institute of Technology, Cambridge, MA 02139, \texttt{\url{http://math.mit.edu/~levine}}}
\address{David B. Wilson, Microsoft Research, Redmond, WA 98052, \texttt{\url{http://dbwilson.com}}}

\date{March 9, 2010}

\begin{abstract}
  A popular theory of self-organized criticality relates the critical
  behavior of driven dissipative systems to that of systems with
  conservation.  In particular, this theory predicts that the
  stationary density of the abelian sandpile model should be equal to
  the threshold density of the corresponding fixed-energy sandpile.
  This ``density conjecture'' has been proved for the underlying graph
  $\Z$.  We show (by simulation or by proof) that the density
  conjecture is \emph{false\/} when the underlying graph is any of
  $\Z^2$, the complete graph $K_n$, the Cayley tree, the ladder graph,
  the bracelet graph, or the flower graph.  Driven dissipative
  sandpiles continue to evolve even after a constant fraction of the sand
  has been lost at the sink.  These results cast doubt on the validity
  of using fixed-energy sandpiles to explore the critical behavior of
  the abelian sandpile model at stationarity.
\end{abstract}

\maketitle

\section{Introduction}
In a widely cited series of papers
\cite{absorbing1,absorbing2,absorbing3,absorbing4,absorbing5},
Dickman, Mu{\~n}oz, Vespignani and Zapperi (DMVZ) developed a theory of
self-organized criticality as a relationship between driven
dissipative systems and systems with conservation.  This theory
predicts a specific relationship between the classical abelian
sandpile model of Bak, Tang, and Wiesenfeld \cite{BTW}, a driven system
in which particles added at random dissipate across the boundary, and
the corresponding ``fixed-energy sandpile,'' a closed system in which
the total number of particles is conserved.

In this introduction, we briefly define these two models and explain
the conjectured relationship between them in the DMVZ paradigm of
self-organized criticality. In particular, we focus on the prediction that the stationary density of the driven dissipative model equals the threshold density of the fixed-energy sandpile model. In \sref{sec:Z2}, we present data from large-scale simulations which strongly indicate that this conjecture
is false on the two-dimensional square lattice $\Z^2$.  In
the subsequent sections we expand on the results announced in \cite{shortversion} by examining the conjecture on some simpler
families of graphs in which we can provably refute it.

The difference between the stationary and threshold densities on most of these graphs
is fairly small --- typically on the order of $0.01\%$ to $0.2\%$ ---
which explains why many previous simulations did not uncover them.  The exception is the early experiments by Grassberger and Manna \cite{GM}, who clearly identified this discrepancy, at least in dimensions $4$ and higher.  Later studies focused on dimension $2$ and missed the discrepancy.

In some more recent papers such as \cite{BM}, the DMVZ paradigm is explicitly restricted to stochastic models.  
In other recent papers \cite{VD,CVSD} it is claimed to apply both to stochastic and deterministic sandpiles, although these papers focus on stochastic sandpiles, for the reason that deterministic sandpiles are said to belong to a different universality class.

Despite our contrary findings, we believe that the central idea of the DMVZ paradigm is
a good one: even in the deterministic case, the dynamics of a driven dissipative system should in some way reflect the dynamics of the corresponding conservative system.  Our results point to a somewhat different relationship than that
posited in the DMVZ series of papers: the driven dissipative model
exhibits a second-order phase transition at the threshold density of
the conservative model.  We explain this transition in
\sref{sec:bracelet}.

Bak, Tang, and Wiesenfeld~\cite{BTW} introduced the abelian sandpile as
a model of self-organized criticality; for mathematical background,
see~\cite{frank}.  The model begins with a collection of particles on
the vertices of a finite graph.  A vertex having at least as many
particles as its degree \emph{topples\/} by sending one particle along
each incident edge.  A subset of the vertices are distinguished as
sinks: they absorb particles but never topple.  A single time step
consists of adding one particle at a random site, and then performing
topplings until each non-sink vertex has fewer particles than its degree.
The order of topplings does not affect the outcome~\cite{dhar}.  The
set of topplings that occur before the system stabilizes is called an avalanche.

Avalanches can be decomposed into a sequence of ``waves'' so that each site topples at most once during each wave.  Over time, sandpiles evolve toward a stationary state in which 
the waves exhibit power-law statistics \cite{KLGP}
 (though the full avalanches seem to exhibit multifractal behavior \cite{MST,KMS}).
Power law statistics are a hallmark of criticality, and since the
stationary state is reached apparently without tuning of a parameter,
the model is said to be \emph{self-organized critical}.

To explain how the sandpile model self-organizes to reach the critical
state, Dickman \textit{et al.}~\cite{absorbing1, absorbing3}
introduced an argument which soon became widely accepted: see, for
example, \cite[Ch.~15.4.5]{sornette} and \cite{quant,feyredig,RS}.
Despite the apparent lack of a free parameter, they argued, the
dynamics implicitly involve the tuning of a parameter to a value where
a phase transition takes place. The phase transition is between an
active state, where topplings take place, and a quiescent
``absorbing'' state where topplings have died out.  The parameter is the \emph{density}, the average
number of particles per site.  When the system is quiescent, addition
of new particles increases the density.  When the system is active,
particles are lost to the sinks via toppling, decreasing the
density. The dynamical rule ``add a particle when all activity has
died out'' ensures that these two density changing mechanisms balance
one another out, driving the system to the threshold of instability.

To explore this idea, DMVZ introduced the \emph{fixed-energy sandpile\/} model (FES),
which involves an explicit free parameter $\zeta$, the density of particles.
On a graph with $N$ vertices, the system starts with $\zeta N$ particles at vertices chosen independently and uniformly at random.
Unlike the driven dissipative sandpile described above, there are no sinks and no addition of particles, so the total number of particles is conserved.  Subsequently the system evolves through toppling of unstable sites. Usually the parallel toppling order is chosen:  at each time step, all unstable sites topple simultaneously.   In the mathematical literature, this system goes by the name of \emph{parallel chip-firing}~\cite{BG,BLS}.  Toppling may persist forever, or it may stop after some finite time.  In the latter case, we say that the system \emph{stabilizes}; in the terminology of DMVZ, it reaches an ``absorbing state.''

A common choice of underlying graph is the $n\times n$ square grid
with periodic boundary conditions.  It is believed, and supported by
simulations \cite{stairs}, that there is a \emph{threshold density}
$\zeta_c$, such that for $\zeta<\zeta_c$, the system stabilizes with
probability tending to~$1$ as $n \to \infty$; and for $\zeta>\zeta_c$,
with probability tending to~$1$ the system does not stabilize.

For the driven dissipative sandpile on the $n\times n$ square grid, as
$n \to \infty$ the stationary measure has an infinite-volume limit
\cite{AJ}, which is a measure on sandpiles on $\Z^2$.  It turns out
that one gets the same limiting measure whether the grid has periodic
or open boundary conditions, and whether there is one sink vertex or
the whole boundary serves as a sink \cite{AJ} (see also
\cite{pemantle} for the corresponding result on random spanning
trees).  The statistical properties of this limiting measure have been
much studied \cite{priezzhevheights,JPR}.  Of particular interest is the \emph{stationary density\/} $\zeta_s$ of $\Z^2$, defined as the expected number of particles at a fixed site.  Grassberger conjectured that $\zeta_s$ is exactly $17/8$, and it is now known that $\zeta_s = 17/8 \pm 10^{-12}$ \cite{JPR}.  

DMVZ believed that the combination of driving and dissipation in the
classical abelian sandpile model should push it toward the critical
density $\zeta_c$ of the fixed-energy sandpile.  This leads to a
specific testable prediction, which we call the Density Conjecture.

\begin{conjecture}[Density Conjecture, \cite{absorbing4}] \label{conj:density}
On the square grid, $\zeta_c = 17/8$.
\end{conjecture}

One can also formulate a density conjecture for more general graphs,
where it takes the form $\zeta_c = \zeta_s$.  We give precise
definitions of the densities $\zeta_c$ and $\zeta_s$ in
section~\ref{sec:Z2}.

Vespignani \textit{et al.}~\cite{absorbing4} write of the fixed-energy
sandpile on the square grid, ``the system turns out to be critical
only for a particular value of the energy density equal to that of the
stationary, slowly driven sandpile.''  They add that the threshold
density $\zeta_c$ of the fixed-energy sandpile is ``the only possible
stationary value for the energy density'' of the driven dissipative
model.  In simulations they find $\zeta_c = 2.1250(5)$, adding in a
footnote ``It is likely that, in fact, 17/8 is the exact result.''

Mu\~{n}oz \textit{et al.}~\cite{absorbing5} have also expressed this view,
asserting that ``FES are found to be critical only for a particular
value $\zeta=\zeta_c$ (which as we will show turns out to be identical
to the stationary energy density of its driven dissipative
counterpart).''

\begin{table}[t!]
\newcommand{\stddev}{0.0000004}
\begin{center}
\parbox{\textwidth}{
\hfill
\begin{tabular}[b]{rcc}
 $n$    &   \#trials & estimate of $\zeta_c(\Z_n^2)$ \\
\hline
   $64$ &   $2^{28}$ & $2.1249561 \pm \stddev$\\
  $128$ &   $2^{26}$ & $2.1251851 \pm \stddev$\\
  $256$ &   $2^{24}$ & $2.1252572 \pm \stddev$\\
  $512$ &   $2^{22}$ & $2.1252786 \pm \stddev$\\
 $1024$ &   $2^{20}$ & $2.1252853 \pm \stddev$\\
 $2048$ &   $2^{18}$ & $2.1252876 \pm \stddev$\\
 $4096$ &   $2^{16}$ & $2.1252877 \pm \stddev$\\
 $8192$ &   $2^{14}$ & $2.1252880 \pm \stddev$\\
$16384$ &   $2^{12}$ & $2.1252877 \pm \stddev$\\
\\
\end{tabular}
\hfill
\psfrag{n}[cr][cr]{$n$}
\psfrag{zc(Zn)}[bl][cl]{$\zeta_c(\Z_n^2)$}
\psfrag{2.125288}[Bl][Bl]{$2.125288$}
\psfrag{2.125}[Bl][Bl]{$2.125$}
\includegraphics[width=0.35\textwidth,height=0.32\textwidth]{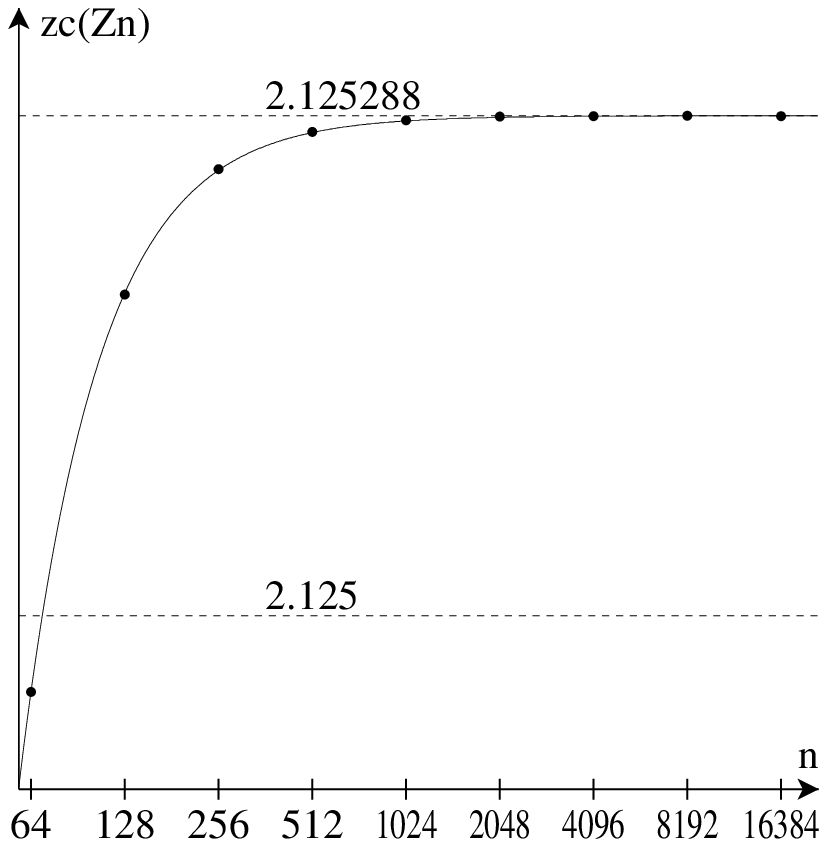}
\hfill\hfill
}
\end{center}
\caption{
  Summary of our fixed-energy sandpile simulations on $n\times n$ tori~$\Z_n^2$, giving our empirical estimate of the threshold density~$\zeta_c(\Z_n^2)$.  The standard deviation in each of our estimates of $\zeta_c(\Z_n^2)$ is $4\times 10^{-7}$.  The data are well approximated by $\zeta_c(\Z_n^2) = 2.1252881\pm 3\times 10^{-7} - (0.390\pm 0.001) n^{-1.7}$, as shown in the graph.  (The error bars are too small to be visible, so the data are shown as points.)  We conclude that the asymptotic threshold density $\zeta_c(\Z^2)$ is $2.125288$ to six decimal places.  In contrast, the stationary density $\zeta_s(\Z^2)$ is $2.125000000000$ to twelve decimal places.
}
\label{table:Z2-density}
\end{table}

Our goal in the present paper is to demonstrate that the density
conjecture is more problematic than it first appears.
Table~\ref{table:Z2-density} presents data from large-scale simulations which
strongly suggest that $\zeta_c$ is close to but not exactly equal to
$17/8$ (see also Table~\ref{table:Z2}).  We further consider several other families of graphs,
including some for which we can determine the exact values $\zeta_c$
and $\zeta_s$ analytically.  We find that they are close, but not
equal.

\begin{table}[h]
\begin{center}
\begin{tabular}{c|c|c}
\hline
graph & $\zeta_s$ & $\zeta_c$ \\
\hline
$\Z$ & {\bf 1} & {\bf 1} \\
$\Z^2$  & $\mathbf{17/8}=2.125$ & $2.125288\ldots$ \\
bracelet & $\mathbf{5/2} = 2.5$ &  $\mathbf{2.496608\ldots}$ \\
flower graph & $\mathbf{5/3} = 1.666667\ldots$ & $\mathbf{1.668898\ldots}$ \\
ladder graph & $\mathbf{\frac74 - \frac{\sqrt{3}}{12}} = 1.605662\ldots$ & $1.6082\ldots$ \\
complete graph  & $n/2 + O(\sqrt{n})$ & $n-O(\sqrt{n \log n})$ \\
3-regular tree  & $\mathbf{3/2}$ & 1.50000\dots \\
4-regular tree  & $\mathbf{2}$ & 2.00041\dots \\
5-regular tree  & $\mathbf{5/2}$ & $2.51167\dots$ \\
\hline
\end{tabular}
\end{center}
\caption{Stationary and threshold densities for different graphs.  Exact values are in bold.
}
\label{table:summary}
\end{table}

All now known information on the threshold density $\zeta_c$ and
stationary density $\zeta_s$ is summarized in
Table~\ref{table:summary}.  The only graph on which the two densities
are known to be equal is $\Z$ \cite{quant,feyredig,FMR}.  On all other
graphs we examined, with the possible exception of the $3$-regular
tree, it appears that $\zeta_c \neq \zeta_s$.

Taken together, these results show that the conclusions of
\cite{absorbing5} that ``FES are shown to exhibit an absorbing state
transition with critical properties coinciding with those of the
corresponding sandpile model'' deserve to be re-evaluated.
One hope of the DMVZ paradigm was that critical features of the driven
dissipative model, such as the exponents governing the distribution of
avalanche sizes and decay of correlations, might be more easily
studied in FES by examining the scaling behavior of these observables
as $\zeta \uparrow \zeta_c$.  However, the failure of the density conjecture
suggests that the two models may not share the same critical features.

As Grassberger and Manna observed~\cite{GM}, the value of the FES threshold density depends on the choice of initial condition.  One might consider a more general version of FES, namely adding $(\zeta-\zeta_0) N$ particles at random to a ``background'' configuration $\tau$ of density $\zeta_0$ already present on the grid.  For example, taking $\tau$ to be the deterministic configuration on $\Z^d$ of $2d-2$ particles everywhere, 
by \cite[Prop~1.4]{FLP} we obtain a threshold density of $\zeta_c = 2d-2$.  Many interesting questions present themselves: for instance, for which background does $\zeta_c$ take the smallest value, and for which backgrounds do we obtain $\zeta_c = \zeta_s$?  It would also be interesting to replicate the phase transition for driven sandpiles (see~\sref{sec:bracelet}) for different background configurations.

\section{Sandpiles on the square grid \texorpdfstring{$\Z^2$}{Z\texttwosuperior}}
\label{sec:Z2}

In this section we give precise definitions of the stationary and threshold densities, and present the results of large-scale simulations on $\Z^2$.  The definitions in this section apply to general graphs, but we defer the discussion of results about other graphs to subsequent sections.

\subsection{The driven dissipative sandpile and the stationary density\texorpdfstring{ $\zeta_s$}{}}

Let $\hat{G} = (V,E)$ be a finite graph, which may have loops and multiple edges.
Let $S \subset V$ be a nonempty set of vertices, which we will call \emph{sinks}.  The presence of sinks distinguishes the driven dissipative sandpile from its fixed-energy counterpart.  To highlight this distinction, throughout the paper, graphs denoted with a ``hat'' as in~$\hat{G}$ have sinks, and those without a hat as in~$G$ do not.

For vertices $v,w\in V$, write $a_{v,w}=a_{w,v}$ for the number of edges connecting~$v$ and~$w$, and	
	\[ d_v = \sum_{w \in V} a_{v,w}. \]
for the number of edges incident to~$v$.  A \emph{sandpile\/} (or ``configuration'') $\eta$ on $\hat{G}$ is a map
 	\[ \eta : V \to \Z_{\geq 0}. \]
We interpret $\eta(v)$ as the number of sand particles at the vertex~$v$; we will sometimes call this number the \emph{height\/} of~$v$ in $\eta$.

 A vertex~$v \notin S$ is called \emph{unstable\/} if $\eta(v) \geq d_v$.  An unstable vertex can \emph{topple\/} by sending one particle along each edge incident to $v$.  Thus, toppling~$v$ results in a new sandpile $\eta'$ given by
	\[ \eta' = \eta + \Delta_v \]
where
	\[ \Delta_v(w) = \begin{cases} a_{v,w}, & v\neq w \\
							a_{v,v} - d_v, &v=w. \end{cases} \]
Sinks by definition are always stable, and never topple.  If all vertices are stable, we say that ~$\eta$ is stable.

Note that toppling a vertex may cause some of its neighbors to become unstable.  The \emph{stabilization\/} $\eta^\circ$ of $\eta$ is a sandpile resulting from toppling unstable vertices in sequence, until all vertices are stable.  By the \emph{abelian property\/}~\cite{dhar}, the stabilization is unique: it does not depend on the toppling sequence.  Moreover, the number of times a given vertex topples does not depend on the toppling sequence.

\begin{figure}
\centering
\includegraphics[width=0.3\textwidth]{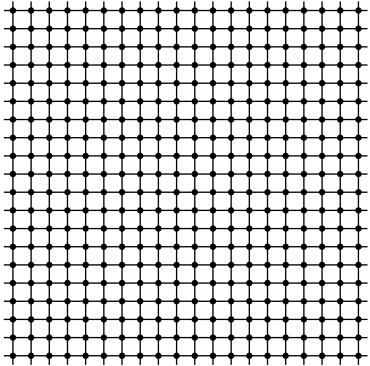}
\caption{The square grid $\Z^2$.}
\label{fig:Z2}
\end{figure}

The most commonly studied example is the $n\times n$ square grid graph, with the boundary sites serving as sinks (\fref{fig:Z2}).

The driven dissipative sandpile model is a continuous time Markov chain $(\eta_t)_{t\geq 0}$ whose state space is the set of stable sandpiles on $\hat{G}$.  Let $V' = V\setminus S$ be the set of vertices that are not sinks.  At each site $v \in V'$, particles are added at rate $1$.  When a particle is added, topplings occur instantaneously to stabilize the sandpile.
Writing $\sigma_t(v)$ for the total number of particles added at~$v$ before time $t$, we have by the abelian property
\[ \eta_t = (\sigma_t)^\circ. \] Note that for fixed~$t$, the random
variables $\sigma_t(v)$ for $v \in V'$ are independent and have the
Poisson distribution with mean~$t$.

The model just described is most commonly known as the abelian
sandpile model (ASM), but we prefer the term ``driven dissipative'' to
distinguish it from the fixed-energy sandpile described below, which
is also a form of ASM.  ``Driven'' refers to the addition of
particles, and ``dissipative'' to the loss of particles absorbed by
the sinks.

Dhar~\cite{dhar} developed the \emph{burning algorithm\/} to
characterize the recurrent sandpile states, that is, those sandpiles
$\eta$ for which, regardless of the initial state, \[ \Pr(\eta_t = \eta \mbox{ for some } t)=1. \]

\begin{lemma}[Burning Algorithm \cite{dhar}]
  A sandpile $\eta$ is recurrent if and only if every non-sink vertex
  topples exactly once during the stabilization of $\eta+\sum_{s\in S}
  \Delta_s$, where the sum is over sink vertices $S$.
\end{lemma}

The recurrent states form an abelian group under the operation of addition followed by stabilization.  In particular, the stationary distribution of the Markov chain $\eta_t$ is uniform on the set of recurrent states.

The combination of driving and dissipation organizes the system into a critical state.  To measure the density of particles in this state, we define the \emph{stationary density\/} $\zeta_s(\hat{G})$ as
	\[
	\zeta_s (\hat{G}) = \EE_{\mu} \left[ \frac{1}{\#V'} \sum_{v \in V'} \eta(v)  \right]
	\]
where $V' = V\setminus S$, and $\mu$ is the uniform measure on recurrent sandpiles on $\hat{G}$.  The stationary density has another expression in terms of the Tutte polynomial of the graph obtained from~$\hat{G}$ by collapsing the set~$S$ of sinks to a single vertex; see \sref{sec:completegraph}.

Most of the graphs we will study arise naturally as finite subsets of a locally finite graph $\Gamma$, i.e., $\Gamma$ is a countably infinite graph in which every vertex has finite degree.  (We also study the complete graph and the flower graph, which do not arise in this way.)  Let $\hat{G}_n$ for $n\geq 1$ be a nested family of finite induced subgraphs with $\bigcup \hat{G}_n = \Gamma$.  As sinks in $\hat{G}_n$ we take the set of boundary vertices
 	\[ S_n = \hat{G}_n - \hat{G}_{n-1}. \]
In cases where the free and wired limits are different, such as on regular trees, we will choose a sequence $\hat{G}_n$ corresponding to the wired limit.
We denote by~$\mu_n$ the uniform measure on recurrent configurations on~$\hat{G}_n$.

We are interested in the \textit{stationary density}
\[ \zeta_s (\Gamma) := \lim_{n \to \infty} \zeta_s( \hat{G}_n).
\]
When $\Gamma=\Z^d$,
it is known that the infinite-volume limit of measures $\mu = \lim_{n \to \infty} \mu_n$ exists and is translation-invariant \cite{AJ}.
In this case it follows that the limit defining $\zeta_s(\Gamma)$ exists and equals
\[
\zeta_s = \E_\mu[\eta(0)].
\]
For other families of graphs we consider, we will show that the limit defining $\zeta_s(\Gamma)$ exists.

Much is known about the limiting measure $\mu$ in the case
$\Gamma=\Z^2$.  The following expressions have been obtained for
$\zeta_s$ and the single site height probabilities.  The symbol $\ceq$ denotes expressions that are rigorous up to a conjecture \cite{JPR} that a
certain integral, numerically evaluated as $0.5 \pm10^{-12}$, is
exactly $1/2$.
\begin{align*}
 \zeta_s(\Z^2) &\ceq 17/8 \quad \text{\cite{JPR}},\\
 \mu\{\eta(x)=0\}&=\textstyle\frac{2}{\pi^2} - \frac{4}{\pi^3} \quad \text{\cite{MD}}, \\
 \mu\{\eta(x)=1\}&\ceq\textstyle\frac14-\frac{1}{2\pi}-\frac{2}{\pi^2}+\frac{12}{\pi^3} \quad \text{\cite{priezzhevheights,JPR}}, \\
 \mu\{\eta(x)=2\}&\ceq\textstyle\frac38+\frac1\pi-\frac{12}{\pi^3} \quad \text{\cite{JPR}, and} \\
 \mu\{\eta(x)=3\}&\ceq\textstyle\frac38-\frac1{2\pi}+\frac1{\pi^2}+\frac4{\pi^3} \quad \text{\cite{JPR}}.
\end{align*}

\subsection{The fixed-energy sandpile model and the threshold density\texorpdfstring{ $\zeta_c$}{}}

Next we describe the fixed-energy sandpile model, in which the driving and dissipation are absent, and the total number of particles is conserved.  As before, let $G$ be a finite graph, possibly with loops and multiple edges.  Unlike the driven dissipative model, we do not single out any vertices as sinks.  The fixed-energy sandpile evolves in discrete time: at each time step, all unstable vertices topple once in parallel.  Thus the configuration $\eta_{j+1}$ at time $j+1$ is given by
	\[ \eta_{j+1} = \eta_j + \sum_{v \in U_j} \Delta_v \]
where
	\[ U_j = \{v \in V \,:\, \eta_j(v) \geq d_v\} \]
is the set of vertices that are unstable at time~$j$.  We say that $\eta_0$ \emph{stabilizes\/} if toppling eventually stops, i.e. $U_j = \varnothing$ for all sufficiently large~$j$.

If $\eta_0$ stabilizes, then there is some site that never topples
\cite{tardos} (see also \cite[Theorem 2.8, item 4]{FMR} and
\cite[Lemma 2.2]{FLP}).  Otherwise, for each site $x$, let $j(x)$ be
the last time $x$ topples.  Choose a site $x$ minimizing $j(x)$.  Then
each neighbor $y$ of $x$ has $j(y) \geq j(x)$, so $y$ topples at least
once at or after time $j(x)$.  Thus $x$ receives at least $d_x$
additional particles and must topple again after time $j(x)$, a
contradiction.  This criterion is very useful in simulations: as soon
as every site has toppled at least once, we know that the system will
not stabilize.

Let $(\sigma_\lambda(v))_{\lambda \geq 0}$ be a collection of independent Poisson point processes of intensity~$1$, indexed by the vertices of $G$.  So each $\sigma_\lambda(v)$ has the Poisson distribution with mean $\lambda$.  We define the \emph{threshold density} of $G$ as
	\[ \zeta_c(G) = \EE \Lambda_c, \]
where
	\[ \Lambda_c = \sup \{\lambda \,:\, \sigma_\lambda \mbox{ stabilizes} \}.  \]
We expect that $\Lambda_c$ is tightly concentrated around its mean when $G$ is large.
Indeed, if $\Gamma$ is an infinite vertex-transitive graph, then the event that~$\sigma_\lambda$ stabilizes on $\Gamma$ is translation-invariant.  By the ergodicity of the Poisson product measure, this event has probability~$0$ or~$1$.  Since this probability is monotone in $\lambda$, there is a (deterministic) threshold density~$\zeta_c(\Gamma)$, such that
	\[ \Pr [\sigma_\lambda \mbox{ stabilizes on } \Gamma] =
		\begin{cases} 1, & \lambda < \zeta_c(\Gamma) \\
				       0, & \lambda > \zeta_c(\Gamma).
		\end{cases} \]
We expect the threshold densities on natural families of finite graphs to satisfy a law of large numbers such as the following.

\begin{conjecture}
  With probability~$1$,
    \[ \Lambda_c(\Z_n^2) \to \zeta_c(\Z^2) \quad \mbox{\em as } n \to \infty. \]
\end{conjecture}

Previous simulations ($n = 160$ \cite{absorbing1}; $n=1280$ \cite{absorbing2}) to estimate the threshold density $\zeta_c(\Z^2)$ found a value of $2.125$, in agreement with the stationary density $\zeta_s (\Z^2)\ceq 17/8$.    By performing larger-scale simulations, however, we find that $\zeta_c$ exceeds $\zeta_s$.

Table~\ref{table:Z2} summarizes the results of our simulations, which
indicate that $\zeta_c(\Z^2)$ equals $2.125288$ to six decimal places.  In each
random trial, we add particles one at a time at uniformly
random sites of the $n\times n$ torus.  After each addition, we
perform topplings until either all sites are stable, or every site has toppled at least once since the last addition.  For deterministic sandpiles on a
connected graph, if every site topples at least once, the system will
never stabilize \cite{FMR,FLP}.  
We record $m/n^2$ as an empirical estimate of the
threshold density, where~$m$ is the maximum number of particles for
which the configuration stabilizes.  We then average these empirical
estimates were over many independent trials.  The one-site marginals
we report are obtained from the stable configuration just before the
$(m+1)^{\text{st}}$ particle was added, and the number of topplings
reported is the total number of topplings required to stabilize the
first $m$ particles.

We used a random number generator based on the Advanced Encryption
Standard (AES-256), which has been found to exhibit good randomness
properties.  Our simulations were conducted on a High Performance
Computing (HPC) cluster of computers.

\begin{table}
\begin{center}
\small
\begin{tabular}{|r|r|c|c|c|c|c|c|}
\hline
grid size\ & \multirow{2}{*}{\#samples} & \multirow{2}{*}{$\zeta_c(\Z_n^2)$} &
 \multicolumn{4}{c|}{distribution of height $h$ of sand} & \#topplings \\
\cline{4-7}
 ($n^2$) & & & $\Pr[h=0]$ & $\Pr[h=1]$ & $\Pr[h=2]$ & $\Pr[h=3]$ &  $\div n^3$ \\
\hline
   $64^2$ & 268435456 & 2.124956 & 0.073555 & 0.173966 & 0.306447 & 0.446032 & 0.197110 \\
  $128^2$ &  67108864 & 2.125185 & 0.073505 & 0.173866 & 0.306567 & 0.446062 & 0.197808 \\
  $256^2$ &  16777216 & 2.125257 & 0.073488 & 0.173835 & 0.306609 & 0.446068 & 0.198789 \\
  $512^2$ &   4194304 & 2.125279 & 0.073481 & 0.173826 & 0.306626 & 0.446067 & 0.200162 \\
 $1024^2$ &   1048576 & 2.125285 & 0.073479 & 0.173822 & 0.306633 & 0.446066 & 0.201745 \\
 $2048^2$ &    262144 & 2.125288 & 0.073478 & 0.173821 & 0.306635 & 0.446065 & 0.203378 \\
 $4096^2$ &     65536 & 2.125288 & 0.073477 & 0.173821 & 0.306637 & 0.446064 & 0.205323 \\
 $8192^2$ &     16384 & 2.125288 & 0.073477 & 0.173821 & 0.306638 & 0.446064 & 0.206475 \\
$16384^2$ &      4096 & 2.125288 & 0.073478 & 0.173821 & 0.306638 & 0.446064 & 0.208079 \\
\hline
\multicolumn{2}{|c|}{$\Z^2$ (stationary)}
                      & 2.125000 & 0.073636 & 0.173900 & 0.306291 & 0.446172 & \multicolumn{1}{c}{} \\
\cline{1-7}
\end{tabular}
\end{center}
\caption{
  Fixed-energy sandpile simulations on $n\times n$ tori $\Z_n^2$.  The third column gives our empirical estimate of the threshold density~$\zeta_c(\Z_n^2)$.  The next four columns give the empirical distribution of the height of a fixed vertex in the stabilization $(\sigma_\lambda)^\circ$, for $\lambda$ just below $\Lambda_c$.  Each estimate of the expectation $\zeta_c(\Z_n^2)$ and of the marginals $\Pr[h=i]$ has standard deviation less than $4 \cdot 10^{-7}$.
  The total number of topplings needed to stabilize $\sigma_\lambda$ appears to scale as $n^3$.  
}
\label{table:Z2}
\end{table}

\section{Sandpiles on the bracelet}
\label{sec:bracelet}

Next we examine a family of graphs for which we can determine $\zeta_c$ and $\zeta_s$ exactly and prove that they are not equal.  Despite this inequality, we show that an interesting connection remains between the driven dissipative and conservative dynamics: the threshold density of the conservative model is the point at which the driven dissipative model begins to lose a macroscopic amount of sand to the sink.

\begin{figure}[b]
\centering
\includegraphics{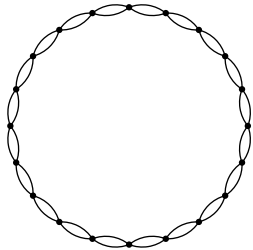}
\caption{The bracelet graph $B_{20}$.}
\label{fig:bracelet}
\end{figure}

The \textit{bracelet graph\/} $B_n$ (\fref{fig:bracelet}) is a multigraph with vertex set
$\Z_n$ (the $n$-cycle) with the usual edge set $\{(i,i+1\bmod n) \,:\, 0 \leq
i <n\}$ doubled.  Thus all vertices have degree~$4$.
The graph~$\hat B_n$ is the same, except that vertex $0$ is distinguished
as a sink from which particles disappear from the system.
We denote by $B_\infty$ the infinite path $\Z$ with doubled edges.

For $\lambda>0$, let $\sigma_\lambda$ be the configuration with
Poisson($\lambda$) particles independently on each site of~$\hat B_n$.  Let $\eta_\lambda = (\sigma_\lambda)^\circ$ be the stabilization of $\sigma_\lambda$, and let
	\[ \rho_n(\lambda) = \frac{1}{n-1} \sum_{x=1}^{n-1} \eta_\lambda(x) \]
be the final density.  The following theorem gives the threshold and stationary densities of the infinite bracelet graph~$B_\infty$, and identifies the $n\to \infty$ limit of the final density $\rho_n(\lambda)$ as a function of the initial density $\lambda$.

\begin{theorem}
\label{braceletmain}
For the bracelet graph,
\begin{enumerate}
\item The threshold density $\zeta_c(B_\infty)$
is the unique positive root of $\zeta = \frac52 - \frac12 e^{-2\zeta}$
(numerically, $\zeta_c = 2.496608$).

\item The stationary density $\zeta_s(B_\infty)$ is $5/2$.

\item
$\rho_n(\lambda) \to \rho(\lambda)$ in probability as $n \to \infty$, where
        \[ \rho(\lambda) = \min\left (\lambda, \frac{5-e^{-2\lambda}}{2}\right) =   \begin{cases} \lambda, & \lambda \leq \zeta_c \\ \frac{5-e^{-2\lambda}}{2}, & \lambda>\zeta_c. \end{cases} \]
\end{enumerate}
\end{theorem}

\begin{figure}
\centering
\psfrag{r}{$\rho$}
\psfrag{l}{$\lambda$}
\psfrag{rc}{$\zeta_c$}
\psfrag{rs=5/2}{$\zeta_s=5/2$}
\psfrag{curve}{$\frac{5-e^{-2\lambda}}{2}$}
\includegraphics[width=\textwidth]{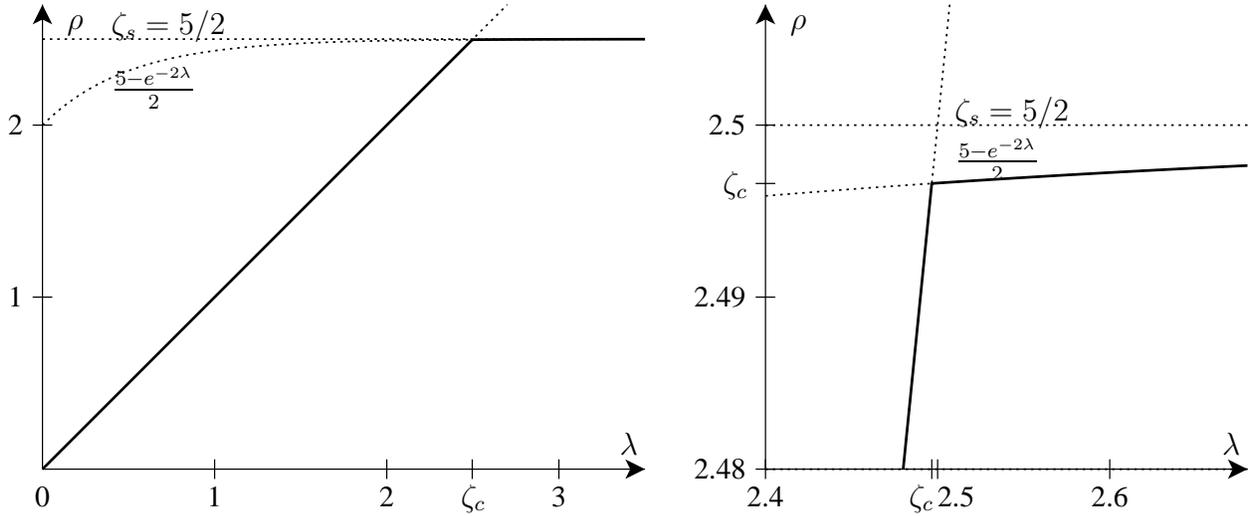}
\caption{Density $\rho(\lambda)$ of the final stable configuration as a function of initial density $\lambda$, for the driven sandpile on the bracelet graph $\hat B_n$ as $n \to \infty$.  A second-order phase transition occurs at $\lambda=\zeta_c$.  Beyond this transition, the density of the driven sandpile continues to increase, approaching the stationary density $\zeta_s$ from below.}
\label{fig:bracelet-density}
\end{figure}

Part 3 of this theorem shows that the final density
undergoes a second-order phase transition at $\zeta_c$: the derivative of $\rho(\lambda)$ is discontinuous at $\lambda=\zeta_c$ (\fref{fig:bracelet-density}).  Thus
in spite of the fact that $\zeta_s \neq \zeta_c$, there remains a
connection between the conservative dynamics used to define $\zeta_c$
and the driven-dissipative dynamics used to define $\zeta_s$.  For
$\lambda<\zeta_c$, very little dissipation takes place, so the final
density equals the initial density $\lambda$; for $\lambda > \zeta_c$
a substantial amount of dissipation takes place, many particles are
lost to the sink, and the final density is strictly less than the
initial density.  The sandpile continues to evolve as $\lambda$ increases
beyond~$\zeta_c$; in particular its density keeps changing.

We believe that this phenomenon is widespread.  As evidence, in
\sref{sec:flower} we introduce the ``flower graph,'' which looks
very different from the bracelet, and prove (in \tref{flowermain}) that a
similar phase transition takes place there.

For the proof of \tref{braceletmain}, we compare the dynamics of pairs
of particles on the bracelet graph to single particles on $\Z$.  At
each vertex $x$ of the bracelet, we group the particles starting
at~$x$ into pairs, with one ``passive'' particle left over
if $\sigma_\lambda(x)$ is odd.  Since all edges in the bracelet are
doubled, we can ensure that in each toppling the two particles
comprising a pair always move to the same neighbor, and that
the passive particles never move.  The toppling dynamics of the pairs
are equivalent to the usual abelian sandpile dynamics on~$\Z$.

We recall the relevant facts about one-dimensional sandpile dynamics:
\begin{itemize}
\item In any recurrent configuration
on a finite interval of $\Z$, every site has height~$1$, except for at most one site of height~$0$. Therefore,~$\zeta_s = 1$ \cite{frank}.
\item On
$\Z$, an initial configuration distributed according to a nontrivial product measure with mean $\lambda$
stabilizes almost surely (every site topples only finitely many times) if $\lambda < 1$, while it almost surely does not stabilize (every site topples infinitely often) if $\lambda \geq 1$ \cite{FMR}. Thus, $\zeta_c = 1$.
 \end{itemize}
 
\begin{proof}[Proof of \tref{braceletmain} parts 1 and 2.]
  Given $\lambda>0$, let $\lambda^*$ be the pair density $\EE
  \lfloor\sigma_\lambda(x)/2\rfloor$, and let
   \[ p_\odd(\lambda) = e^{-\lambda}\sum_{m \geq 0} \frac{\lambda^{2m+1}}{(2m+1)!} =  \frac 12 (1-e^{-2\lambda}). \]
  be the probability that a Poisson($\lambda$) random variable is odd.  Then $\lambda$ and $\lambda^*$ are related by
   \begin{equation} \label{eq:lambdastar}
   \lambda = 2\lambda^* + p_\odd(\lambda).
   \end{equation}
The configuration $\sigma_\lambda$ stabilizes on $B_\infty$ if and only if the pair configuration $\sigma_\lambda^*$ stabilizes on~$\Z$.  Thus $\zeta_c(B_\infty)^* = \zeta_c(\Z)$.  Setting $\lambda = \zeta_c(B_\infty)$ in (\ref{eq:lambdastar}), using the fact that $\zeta_c(\Z)=1$, and that~$\lambda^*$ is an increasing function of $\lambda>0$, we conclude that $\zeta_c(B_\infty)$ is the unique positive root of
	\[ \zeta = 2 + p_\odd(\zeta), \]
or $\zeta = \frac52 - \frac12 e^{-2\zeta}$.
This proves part 1.  

For part~2, by the burning algorithm, a configuration $\sigma$ on $\hat{B}_n$ is recurrent if and only if it has at most one site with fewer than two particles.  Thus, in the uniform measure on recurrent configurations on $\hat{B}_n$,
	\[ \Pr(\sigma(x) = 2) = \Pr(\sigma(x)=3) = \frac12 -\frac{1}{2n}, \qquad \Pr(\sigma(x)=0) = \Pr(\sigma(x)=1) = \frac{1}{2n}. \] 
We conclude that $\zeta_s(\hat{B}_n) = \EE \sigma(x) = \frac52 - \frac{2}{n} \to \frac 52$ as $n \to \infty$. 
\end{proof}

To prove part 3 of \tref{braceletmain}, we use the following lemma, whose proof is deferred to the end of this section.  Let $\hat\Z_n$ be the $n$-cycle with vertex $0$ distinguished as a sink.
Let~$\sigma'_\lambda$ be a sandpile on $\hat\Z_n$ distributed according to a product measure (not necessarily Poisson) of mean $\lambda$.  Let $\eta'_\lambda$ be the stabilization of $\sigma'_\lambda$, and let $\rho'_n(\lambda) = \frac{1}{n-1} \sum_{x=1}^{n-1} \eta'_\lambda(x)$ be the final density after stabilization.

\begin{lemma} \label{linelemma}
On $\hat \Z_n$, we have $\rho'_n(\lambda) \to \min(\lambda,1)$ in probability.
\end{lemma}

\begin{proof}[Proof of \tref{braceletmain}, part 3]
Let $\eta_\lambda$ be the stabilization of $\sigma_\lambda$ on $\hat B_n$, and let $\eta^*_\lambda$ be the stabilization of $\sigma^*_\lambda = \floor{\sigma_\lambda /2}$ on $\hat\Z_n$. Then
\begin{equation}
\label{finalconf}
\eta_\lambda(x) = 2 \eta^*_\lambda(x) + \omega_\lambda(x)
\end{equation}
where $\omega_\lambda(x) = \sigma_\lambda(x) - 2 \sigma^*_\lambda(x)$ is $1$ or $0$ accordingly as $\sigma_\lambda(x)$ is odd or even.
Let
	\[ \rho_n^*(\lambda) =  \frac{1}{n-1} \sum_{x=1}^{n-1} \eta^*_\lambda(x) \]
be the final density after stabilization of $\sigma^*_\lambda$ on $\hat \Z_n$.
Then
	\[ \rho_n(\lambda) = 2\rho_n^*(\lambda) + \frac{1}{n-1} \sum_{x=1}^{n-1} \omega_\lambda(x). \]
By the weak law of large numbers, $\frac{1}{n-1} \sum_{x=1}^{n-1} \omega_\lambda(x) \to p_\odd(\lambda)$ in probability as $n \to \infty$.
If $\lambda<\zeta_c$, then $\lambda^* < 1$, so by \lref{linelemma}, $\rho_n^*(\lambda) \to \lambda^*$ in probability, and hence
	\[ \rho_n(\lambda) \to 2\lambda^* + p_\odd(\lambda) = \lambda \]
in probability.
If $\lambda \geq \zeta_c$, then $\lambda^* \geq 1$, so by \lref{linelemma}, $\rho_n^*(\lambda)  \to 1$ in probability, hence
	\[ \rho_n(\lambda) \to 2 + p_\odd(\lambda) = \frac{5 - e^{-2\lambda}}{2} \]
in probability.  This proves part~3.  
\end{proof}

\begin{proof}[Proof of \lref{linelemma}]
  We may view $\hat\Z_n$ as the path in $\Z$ from
  $a_n=-\floor{n/2}$ to $b_n=\ceiling{n/2}$, with both endpoints
  serving as sinks.
  For $x\in \hat\Z_n$, let $u_n(x)$ be the number of times that $x$ topples during stabilization of the configuration $\sigma'_\lambda$ on $\hat\Z_n$.  Let $u_\infty(x)$ be the number of times  ~$x$ topples during stabilization of $\sigma'_\lambda$ on~$\Z$.
  The procedure of ``toppling in nested volumes'' \cite{FMR} shows that $u_n(x) \uparrow u_\infty(x)$ as $n \to \infty$.

We consider first $\lambda<1$.  In this case $u_\infty(x)$ is finite almost surely (a.s.).  The total number of particles lost to the sinks on $\hat \Z_n$ is $u_n(a_n+1)+u_n(b_n-1)$, so the final density is given by
  	\[ \rho'_n(\lambda) = \frac{1}{n-1} \left[ \sum_{x=a_n+1}^{b_n-1} \sigma_\lambda(x) - u_n(a_n+1) - u_n(b_n-1) \right]. \]
By the law of large numbers, $\frac{1}{n-1} \sum \sigma_\lambda(x) \to \lambda$ in probability as $n \to \infty$.
Since $u_\infty(x)$ is a.s.\ finite, we have $\frac{u_n(a_n+1) + u_n(b_n-1)}{n-1} \to 0$ in probability, so $\rho'_n(\lambda) \to \lambda$ in probability.

Next we consider $\lambda \geq 1$.   In this case we have $u_n(x) \uparrow u_\infty(x) = \infty$, a.s.  
Let $p(n,x)= \Pr(u_n(x)=0)$ be the probability that $x \in \hat \Z_n$ does not topple.  By the abelian property, adding sinks can not increase the number of topplings, so
	\[ p(n,x) \leq p(m, 0) \]
where $m = \min(x-a_n, b_n-x)$.  Let
 	\[ Y_n = \sum_{x=a_n+1}^{b_n-1} 1_{\{u_n(x)=0\}} \] 
be the number of sites in $\hat \Z_n$ that do not topple.  Since $u_n(0) \uparrow \infty$ a.s., we have $p(n,0) \downarrow 0$, hence
	\[ \EE \frac{Y_n}{n} = \frac{1}{n} \sum_{x=a_n+1}^{b_n-1} p(n,x) \leq \frac{2}{n} \sum_{m=1}^{n/2} p(m,0) \to 0 \]
as $n \to \infty$.  Since $Y_n \geq 0$ it follows that $Y_n/n \to 0$ in probability.
		
In an interval where every site toppled, there can be at most one empty site. We have $Y_n+1$ such intervals. Therefore, the number of empty sites is at most $2Y_n+1$. Hence
	\[ \frac{n-2Y_n-2}{n-1}  \leq \rho'_n(\lambda) \leq 1. \]
The left side tends to $1$ in probability, which completes the proof.
\end{proof}

\section{Sandpiles on the complete graph}
\label{sec:completegraph}

Let~$K_n$ be the complete graph
on~$n$ vertices: every pair of distinct vertices is connected by an edge.
In~$\hat K_n$, one vertex is distinguished as the sink.
The maximal stable configuration on~ $\hat K_n$ has density $n-2$, while the minimal recurrent configurations have exactly one vertex of each height $0,1,\ldots,n-2$, hence density~$\frac{n-2}{2}$.  The following result shows that the stationary and threshold densities are quite far apart:~$\zeta_s$ is close to the minimal recurrent density, while~$\zeta_c$ is close to the maximal stable density.

\begin{theorem}
\label{completegraphmain}
	\begin{align*} \zeta_s(\hat K_n) &= \frac{n}{2} + O(\sqrt{n}) \\
			       \zeta_c(K_n) &\geq n - O(\sqrt{n \log n}). \end{align*}
\end{theorem}

The proof uses an expression for the stationary density $\zeta_s$ in terms of the Tutte polynomial, due to C. Merino L\'opez~\cite{MerinoLopez}.  Our application will be to the complete graph, but we state Merino L\'{o}pez' theorem in full generality.
Let $G=(V,E)$ be a connected undirected graph with~$n$ vertices and~$m$ edges.  Let~$v$ be any vertex of~$G$, and write~$\hat G$ for the graph~$G$ with~$v$ distinguished as a sink.  Let~$d$ be the degree of~$v$.

Recall that the Tutte polynomial $T_G(x,y)$
is defined by
	\[ T_G(x,y) = \sum_{A\subseteq E} (x-1)^{c(A)-c(E)} (y-1)^{c(A)+|A|-|V|} \]
where $c(A)$ denotes the number of connected components of the spanning subgraph $(V,A)$.

\begin{theorem}[\cite{MerinoLopez}]
The Tutte polynomial $T_G(x,y)$ evaluated at $x=1$ is given by
	\[ T_G(1,y) = y^{d - m} \sum_{\sigma} y^{|\sigma|} \]
where the sum is over all recurrent sandpile configurations $\sigma$ on $\hat G$, and $|\sigma|$ denotes the number of particles in $\sigma$.
\end{theorem}

Differentiating and evaluating at $y=1$, we obtain
	\begin{equation}
	\label{eq:partialtutte}
	\left.\frac{d}{dy} T_G(1,y)\right|_{y=1} = \sum_{\sigma} (d-m+|\sigma|).
	\end{equation}
Referring to the definition of the Tutte polynomial, we see that $T_G(1,1)$ is the number of spanning trees of~$G$, and that the left side of~(\ref{eq:partialtutte}) is the number of spanning unicyclic subgraphs of~$G$.
 (In evaluating~$T_G$ at $x=y=1$, we interpret~$0^0$ as~$1$.)
The number of recurrent configurations equals the number of spanning trees of $G$, so the stationary density
$\zeta_s$ may be expressed as 	
	\[ \zeta_s(\hat G) = \frac{1}{n T_G(1,1)} \sum_{\sigma} |\sigma|. \]
Combining these expressions yields the following:

\begin{corollary}
\label{unicycle}
	\[ \zeta_s(\hat G) = \frac1n \left( m-d + \frac{u(G)}{\kappa(G)} \right) \]
where $\kappa(G)$ is the number of spanning trees of $G$, and $u(G)$ is the number of spanning unicyclic subgraphs of~$G$. 
\end{corollary}

Note that $m-d$ is the minimum number of particles in a recurrent configuration, so the ratio $u(G)/\kappa(G)$ can be interpreted as the average number of excess particles in a recurrent configuration.

Everything so far applies to general connected graphs $G$.
The following is specific for the complete graph.

\begin{theorem}[Wright \cite{wright}]
\label{unicyclesofK_n}
The number of spanning unicyclic subgraphs of~$K_n$ is
	\[ u(K_n) = \left( \sqrt{\frac{\pi}{8}}+o(1) \right) n^{n-\frac12}. \]
\end{theorem}

\begin{proof}[Proof of \tref{completegraphmain}]
For $\hat K_n$ we have
	\[ m-d = \frac{n(n-1)}{2} - (n-1) = \frac{(n-2)(n-1)}{2}. \]
From \cref{unicycle}, \tref{unicyclesofK_n}, and Cayley's formula $\kappa(K_n) = n^{n-2}$, we obtain
	\begin{align*} \zeta_s(\hat K_n) &= \frac{1}{n} \left( \frac{(n-2)(n-1)}{2} + \frac{u(K_n)}{\kappa(K_n)} \right) \\ &=
	\frac{n}{2} + \left(\sqrt{\frac{\pi}{8}} + o(1) \right) \sqrt{n}. \end{align*}
On the other hand, if we let
	\[ \lambda = n - 2 \sqrt{n\log n} \]
and start with $\sigma(v)\sim\,$Poisson$(\lambda)$ particles at each vertex~$v$ of~$K_n$, then for all~$v$
	\[ \Pr[\sigma(v) \geq n] < \frac{1}{n^2}. \]
So
	\[ \Pr[\sigma(v) \geq n \mbox{ for some } v] < \frac1n; \]
in other words, with high probability no topplings occur at all.  Thus
	\[ \Pr \left( \Lambda_c(K_n) \geq n - 2\sqrt{n \log n} \right) > 1- \frac1n \]
which completes the proof.
\end{proof}

One might guess that the large gap between~$\zeta_s$ and~$\zeta_c$ is related to the small diameter of~$\hat K_n$: since the sink is adjacent to every vertex, its effect is felt with each and every toppling.  This intuition is misleading, however, as shown by the lollipop graph~$\hat L_n$ consisting of~$K_n$ connected to a path of length~$n$, with the sink at the far end of the path.
Since~$L_n$ has the same number of spanning trees and unicyclic subgraphs as~$K_n$, we have by \cref{unicycle}
	\[ \zeta_s(\hat L_n) = \frac{1}{2n} \left( \overbrace{\frac{n(n-1)}{2} + n}^m-1 + \frac{u(L_n)}{\kappa(L_n)} \right) = \frac{n}{4} + O(\sqrt{n}). \]
On the other hand, by first stabilizing the vertices on the path, close to half of which end up in the sink without reaching the $K_n$, it is easy to see that with high probability
	\[ \Lambda_c(L_n) \geq \frac{2n}{3} - O(\sqrt{n \log n}). \]
	

\section{Sandpiles on the flower graph}
\label{sec:flower}

An interesting feature of parallel chip-firing is that further phase transitions appear above the threshold density~$\zeta_c$.  On a finite graph $G=(V,E)$, since the time evolution is deterministic, the system will eventually reach a periodic orbit: for some positive integer $m$, we have $\eta_{t+m} = \eta_t$ for all sufficiently large~$t$.  The \textit{activity density}, $\rho_a$, measures the proportion of vertices that topple in an average time step:

\begin{equation*}
\rho_a (\lambda) = \EE_\lambda \lim_{t \to \infty} \frac{1}{t} \sum_{s = 0}^{t-1} \frac{1}{\#V} \sum_{x \in V} \one_{\{\eta_s(x) \geq d_x\}}.
\end{equation*}

The expectation $\EE_\lambda$ refers to the initial state $\eta_0$, which we take to be distributed according to the Poisson product measure with mean~$\lambda$.  Note that the limit in the definition of $\rho_a$ can also be expressed as a finite average, due to the eventual periodicity of the dynamics.

Bagnoli et al.~\cite{stairs} observed that $\rho_a$ tends to increase with $\lambda$ in a sequence of flat steps punctuated by sudden jumps.  This ``devil's staircase'' phenomenon is so far explained only on the complete graph~\cite{lionel}: The number of flat stairs increases with~$n$, and in the $n\to \infty$ limit there is a stair at each rational number height $\rho_a=p/q$.

On the cycle $\Z_n$ \cite{circle} there are just two jumps: at $\lambda=1$, the activity density jumps from $0$ to $1/2$, and at $\lambda=2$, from $1/2$ to $1$.  For the $n \times n$ torus, simulations \cite{stairs} indicate a devil's staircase, which is still not completely understood despite much effort~\cite{CDVV}.

In this section we study the ``flower'' graph, which was designed with parallel chip-firing in mind: the idea is that a graph with only short cycles should give rise to short period orbits under the parallel chip-firing dynamics.  We find that there are four activity density jumps (\tref{fivestairs}).  In addition, we determine the stationary and threshold densities of the flower graph, and find a second-order phase transition at $\zeta_c$ (\tref{flowermain}).

\begin{figure}
\centering
\includegraphics[width=0.3\textwidth]{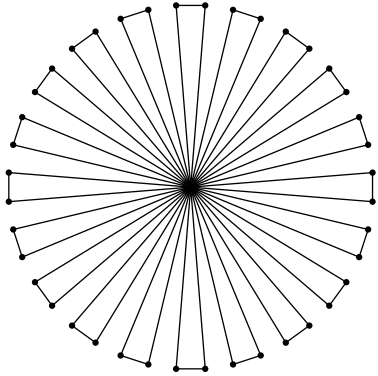}
\caption{The flower graph $F_{20}$.}
\label{fig:flower}
\end{figure}

The flower graph $F_n$ consists of a central site together with $n \geq 1$ petals (\fref{fig:flower}). Each petal consists of two sites connected by an edge, each connected to the central site by an edge.
Thus the central site has degree $2n$, and all other sites have degree 2. The number of sites is $2n+1$.
The graph~$\hat F_n$ is the same, except one petal serves as sink.

Recall that we defined the density of a configuration as the total number of particles, divided by the total number of sites.
Since the flower graph is not regular, the central site has a different expected number of particles than the petal sites.

\begin{prop} \label{period3}
  For parallel chip-firing on the flower graph $F_n$, every configuration has eventual period at most~$3$.
\end{prop}
The proof uses the following two lemmas.
\begin{lemma}[\cite{prisner} Lemma~2.5] \label{maxheight}
  If the eventual period is not $1$, then after some finite time, every site $x$ has height at most $2d_x - 1$.
\end{lemma}
We also use an observation from \cite{stairs}; it was stated and proved there for $\Z^2$, but the same proof works for general graphs.
\begin{lemma}[\cite{stairs}] \label{mirror}
  Let two height configurations $\eta$ and
  $\xi$ be ``mirror images'' of each other, that is, $\eta(x) =
  2d_x -1 -\xi(x)$ for all $x$. Then after performing for each a parallel
  chip-firing time step, the two configurations are again mirror
  images of each other.
\end{lemma}

\begin{proof}[Proof of \pref{period3}]
  Suppose at time $t$ the model has settled into periodic orbit, and
  the period is not $1$.  Then by \lref{maxheight}, at time $t$
  every petal site has height at most $3$.  Say that a petal is in
  state $ij$ if it has $i$ particles at one site and $j$ particles at the
  other site.  \textit{A priori} there are 10 possible petal states,
  listed below, where each one has two possible successor states,
  depending on whether or not the central site is stable. If a petal
  is in state $ij$, then by S($ij$) we denote the state that it is in
  after one time step in which the central site does not topple, and
  likewise by U($ij$) after one time step in which the central site
  topples.
  \begin{center}
  \begin{tabular}{ccc} \small
  state & S(state) & U(state) \\
  \hline
  00 &  00 & 11 \\
  01 &  01 & 12 \\
  02 &  01 & 12 \\
  03 &  11 & 22 \\
  11 &  11 & 22 \\
  12 &  02 & 13 \\
  13 &  12 & 23 \\
  22 &  11 & 22 \\
  23 &  12 & 23 \\
  33 &  22 & 33 \\
  \end{tabular}
  \end{center}
  From this we see that a petal will be in state 00 only if the central
  site is always stable, and consequently each site is always stable,
  in which case the period is 1.  Similarly, petal state 33 only
  occurs if the central site is unstable each step, in which case each
  site must be unstable each step, and the period is again 1.  State
  03 is not a successor of any state of these states, so it will not
  be a periodic petal state either.  Thus the set of allowed periodic
  petal states is $\{01, 02, 11, 12, 13, 22, 23\}$.

  If the central site is stable every other time step, then the possible
  petal states are $12\to02\to12$, $22\to11\to22$, and $13\to12\to13$,
  each of which has period 2. Then the period of the entire configuration is 2.

  Thus if the period is larger than 2, the central site must be stable
  for at least two consecutive time steps, or else unstable for at
  least two consecutive time steps.  We will label a time step S if
  the central site is stable in that time step, otherwise we label it
  U. So, if the period is larger than 2 we will see SS or UU in the
  time evolution. In the latter case, we can study the mirror image,
  which will have the same period, and for which we will see SS.

  Eventually the central site must be unstable again, since otherwise
  the period would be~1. Therefore, we can examine three time steps
  labeled SSU. Examining the evolution of the central site together
  with the petals, we see
  \begin{center}
  \begin{tabular}{ccccccc}
  S &S &U &\\
  01, 02&01&01&12\\
  12, 23&02&01&12\\
  11, 22&11&11&22\\
  23   &12&02&12\\
  \end{tabular}
  \end{center}

  Whenever we have SSU, during the second and third time steps each
  petal contributes at most two particles to the central site, while the
  central site topples, so the central site must again be stable.  Thus
  SSUU cannot occur, and we see SSUS.

  There are two cases for what the central site does next. Let us
  first consider SSUSU.
  \begin{center}
  \begin{tabular}{ccccccc}
  S & S & U & S & U &&\\
  01, 02&01&01&12&02&12&\\
  12, 23&02&01&12&02&12&\\
  11, 22&11&11&22&11&22&\\
  23&12&02&12&02&12&\\
  \end{tabular}
  \end{center}
  During the last two time steps, each petal contributes exactly 2
  particles to the central site, and the central site topples once. Thus
  after two time steps not only the petals, but also the central site
  is in the same state. Therefore, the period becomes 2.

  Next we consider SSUSS.  At this stage each petal is in state 01 or
  11, so if there were yet another S, the sandpile would be periodic
  with period 1. So we see SSUSSU, and because SSUU is forbidden, we
  conclude that we see SSUSSUS.
  \begin{center}
  \begin{tabular}{ccccccc}
  S&S&U&S&S&U&S\\
  01, 02&01&01&12&02&01&12\\
  12, 23&02&01&12&02&01&12\\
  11, 22&11&11&22&11&11&22\\
  23&12&02&12&02&01&12\\
  \end{tabular}
  \end{center}
  At the time of the third S, each petal is in state 12 or 22.  Between
  the third S and the fifth S, each petal contributes exactly two
  particles to the central site and returns to the same state, while the
  central site topples once.  Thus the configuration is periodic with
  period 3.
\end{proof}

We conclude from the above case analysis that the activity $\rho_a$ is always one of $0$, $1/3$, $1/2$, $2/3$, or $1$.
Table~\ref{periodicstates} summarizes the behavior of the periodic sandpile states for different values of $\rho_a$.

\begin{table}[hb]
\begin{center}
\begin{tabular}{|c|c|ccc|cc|ccc|c|}
\cline{2-11}
\multicolumn{1}{c|}{} & \multicolumn{10}{c|}{periodic sandpile states} \\
\cline{1-11}
activity $\rho_a$ &
\multicolumn{1}{|c|}{0} & \multicolumn{3}{c|}{1/3} & \multicolumn{2}{c|}{1/2} & \multicolumn{3}{c|}{2/3} & \multicolumn{1}{c|}{1} \\
\hline
central site &
S & S & S & U & S & U & S & U & U & U \\
\multirow{3}{*}{\begin{rotatebox}{90}{petals}\end{rotatebox}} &
 01& 12& 02& 01& 12& 02& 23& 12& 13& $\geq 22$  \\ &
11& 22& 11& 11& 22& 11& 22& 11& 22&  \\ &
00 &   &   &   & 13& 12&   &   &   &   \\
\hline
\end{tabular}
\end{center}
\caption{Behavior of the central site and petals as a function of the activity $\rho_a$.}
\label{periodicstates}
\end{table}

The following theorem shows that parallel chip-firing on the flower graph exhibits four distinct phase transitions where the activity $\rho_a$ jumps in value:  For each $\alpha \in \{0,\frac13,\frac12,\frac23,1\}$, there is a nonvanishing interval of initial densities $\lambda$ where $\rho_a=\alpha$
asymptotically almost surely.

\begin{theorem}
\label{fivestairs}
  Let $\zeta_c$ be the unique root of $\frac 53 + \frac 13 e^{-3\zeta} =
  \zeta$, and let $\zeta'_c$ be the unique positive root of $\frac{10}{3} - \frac13
  e^{-3\zeta} = \zeta$.  (Numerically, $\zeta_c=1.6688976\dots$ and
  $\zeta'_c=3.3333182\dots$.)  With probability tending to~$1$ as~$n \to \infty$, the activity density~$\rho_a$ of parallel chip-firing on the flower graph~$F_n$ is given by
  \[ \rho_a = \begin{cases}
0,   & \text{\em if  $0 \leq \lambda < \zeta_c$} \\
1/3, & \text{\em if  $\zeta_c < \lambda < 2$} \\
1/2, & \text{\em if  $2 < \lambda < 3$} \\
2/3, & \text{\em if  $3 < \lambda < \zeta'_c$} \\
1   & \text{\em if  $\zeta'_c < \lambda$.}
\end{cases}
\]
\end{theorem}

\begin{proof}
  In a given petal, let $X$ denote the difference modulo $3$ of the number of particles on the two sites of the petal.  Observe that $X$ is unaffected by toppling.  Let
  $Z$ denote the number of petals for which $X=0$, and $R$ denote the
  total number of particles, in a given initial configuration.  Using
  Table~\ref{periodicstates}, we can relate $Z$, $R$, and the
  activity $\rho_a$.

  When $\rho_a=0$, we have less than $2n$ particles at the central site,
  at most two particles for the $Z$ petals of type $X=0$, and exactly one particle
  for the other $n-Z$ petals, so $ 0 \leq R < 2n + 2Z + (n-Z) =
  3n + Z$.

  When $\rho_a=1/3$, by considering the U time step, we have $ R \geq
  2n + 2Z + (n-Z) = 3n + Z$.  By considering the preceding S
  time step, we have $R < 2n + 2Z  + 2(n-Z) = 4n$.

  When $\rho_a=1/2$, by considering the U time step, we have $R \geq
  4n$, and by considering the S time step, we get $R < 2n + 4Z
  + 4(n-Z) = 6n$.

  When $\rho_a=2/3$, by considering the second U step, we have $R \geq
  2n + 4Z + 4(n-Z) = 6n$.  By considering the S time
  step, we have $R < 2n + 4Z + 5(n-Z) = 7n-Z$.

  When $\rho_a=1$, we have $R\geq 2n + 4 Z + 5(n-Z) = 7n - Z$.

  Since for given $n$ and $Z$, these intervals on the values of $R$
  are disjoint, we see that the converse statements hold as well: the
  values of $R$ and $Z$ determine the activity $\rho_a$.  We summarize
  these bounds:
\[ \rho_a = \begin{cases}
0   & \text{if and only if $0 \leq R < 3n+Z$} \\
1/3 & \text{if and only if $3n+Z\leq R<4n$}\\
1/2 & \text{if and only if $4n\leq R<6n$}\\
2/3 & \text{if and only if $6n\leq R<7n-Z$} \\
1   & \text{if and only if $7n-Z\leq R$}.
\end{cases}
\]

  Everything so far holds deterministically; next we use probability to
  estimate $R$ and $Z$.  By the weak law of large numbers,
  $R/n \to 2\lambda$ and $Z/n \to \Pr(X=0)$ in probability.
 Thus, to complete the proof it suffices to show
  	\begin{equation} \label{eq:whatweneeded} \Pr(X=0) = \frac13 (1+2 e^{-3\lambda}). \end{equation}
 
We can think of building the initial configuration $\sigma_\lambda$ by starting with the empty configuration and adding particles in continuous time. Then the value of $X$ for a single petal as particles are added is a continuous time Markov chain on the state space $\{0, \pm 1\}$ with transitions $0 \to \pm 1$ at rate $2$, and $\pm 1 \to 0$ and $\pm 1 \to \pm 1$ each at rate $1$.  Starting in state~$0$, after running this chain for time $\lambda$
we obtain
 	\[ \left[ \begin{array}{c} \Pr(X=0) \\ \Pr(X\neq 0) \end{array} \right] = \exp \left\{ 2 \lambda \left( P-I
	\right) \right\} \left[ \begin{array}{c} 1 \\ 0 \end{array} \right] \] 
where $P= \left[ \begin{smallmatrix} 0 & 1/2 \\ 1 & 1/2 \end{smallmatrix} \right]$, and $I$ is the $2\times 2$ identity matrix.  The eigenvalues of $P-I$ are $0$ and~$-\frac32$, with corresponding eigenvectors $v_1 = \left[ \begin{smallmatrix} 1 \\ 2 \end{smallmatrix} \right]$ and $v_2 = \left[ \begin{smallmatrix} 1 \\ -1 \end{smallmatrix} \right]$. 
Since  $\left[ \begin{smallmatrix} 1 \\ 0 \end{smallmatrix} \right] = \frac13 v_1 + \frac23 v_2$, we obtain~(\ref{eq:whatweneeded}).
\end{proof}

\begin{figure}
\centering
\psfrag{r}{$\rho$}
\psfrag{l}{$\lambda$}
\psfrag{rc}{$\zeta_c$}
\psfrag{rs=5/3}{$\zeta_s=5/3$}
\psfrag{curve}{$\frac{5+e^{-\lambda}}{3}$}
\includegraphics[width=\textwidth]{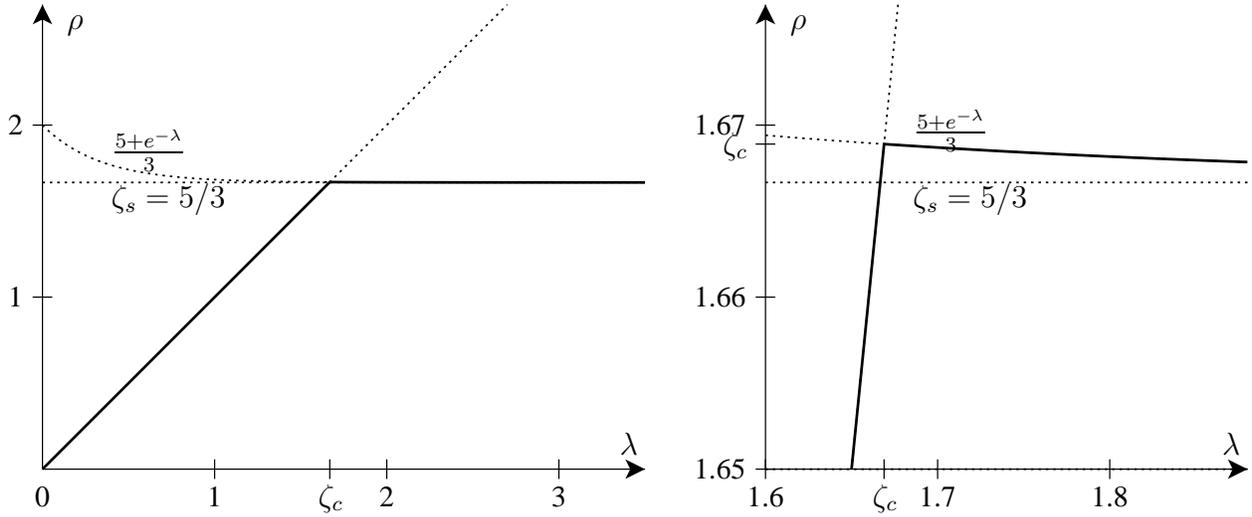}
\caption{Density $\rho(\lambda)$ of the final stable configuration as a function of initial density $\lambda$ on the flower graph $\hat F_n$ for large $n$.  A second-order phase transition occurs at $\lambda=\zeta_c$.  Beyond this transition, the density of the driven sandpile decreases with $\lambda$.}
\label{fig:flower-density}
\end{figure}

The following theorem describes a phase transition in the driven
sandpile dynamics on the flower graph analogous to \tref{braceletmain}
for the bracelet graph.  We remark on one interesting difference
between the two transitions: for $\lambda > \zeta_c$, the final
density $\rho(\lambda)$ is increasing in $\lambda$ for the bracelet,
and decreasing in $\lambda$ for the flower graph.

For $\lambda>0$, let $\sigma_\lambda$ be the configuration with
Poisson($\lambda$) particles independently on each site of~$\hat F_n$.  Let $\eta_\lambda = (\sigma_\lambda)^\circ$ be the stabilization of $\sigma_\lambda$, and let
	\[ \rho_n(\lambda) = \frac{1}{n-1} \sum_{x=1}^{n-1} \eta_\lambda(x) \]
be the final density. 

\begin{theorem} \label{flowermain}
For the flower graph with $n$ petals, in the limit $n\to\infty$ we have
\begin{enumerate}
\item The threshold density $\zeta_c$ is the unique positive root of $\zeta = \frac53 + \frac13 e^{-3\zeta}$.
\item The stationary density $\zeta_s$ is $5/3$.
\item $\rho_n(\lambda) \to \rho(\lambda)$ in probability, where
 \[ \rho(\lambda) = \min\left (\lambda, ~\frac53+\frac13 e^{-\lambda}\right) =   \begin{cases} \lambda, & \lambda \leq \zeta_c \\ \frac53+\frac13 e^{-3\lambda}, & \lambda>\zeta_c. \end{cases} \]
 \end{enumerate}
 \end{theorem}

\begin{proof}
Part 1 follows from \tref{fivestairs}.

For Part 2, we use the burning algorithm.
  In all recurrent
  configurations on $\hat F_n$, the central site has either $2n-1$ or $2n-2$ particles.
  All other sites have at most one particle, and in each petal (except the sink) there is at least one particle.  For each
  petal that is not the sink, there are two possible configurations
  with $1$ particle, and one with $2$ particles.  Each of these occurs
  with equal probability in the stationary state, so the expected
  number of particles in the petals is  $(n-1) (\frac23 \cdot 1 + \frac13 \cdot 2) =
  \frac{4n}{3} + O(1)$ as $n \to \infty$.  Therefore,
  the total density is $\zeta_s = \lim_{n \to \infty} \frac{2n + 4n/3}{2n-1} + o(1) = 5/3$.

For part 3, for the driven dissipative sandpile on $\hat F_n$,
  we first stabilize all the petals, then topple the center site if
  it is unstable, then stabilize all the petals, and so on.  For each
  toppling of the center site, the sandpile loses $O(1)$ particles to the
  sink.  If the center topples at least once, then each petal will be in one
  of the states $11$, $01$, or $10$, after which the number of particles at the center site is $R-n-Z+O(1)$.  Recall from the proof of Theorem~\ref{fivestairs} that $R/n \to 2\lambda$ and $Z/n \to \frac 13 (1+2 e^{-3\lambda})$ in probability.  Thus if $\lambda\leq \zeta_c$, then
   $\frac{R-n+1-Z}{2n} \to \lambda - \frac23 + \frac13 e^{-3\lambda} \leq 1$ in probability, so the sandpile does not lose a macroscopic amount of sand, and $\rho_n(\lambda) \to \lambda$ in probability.

If $\lambda>\zeta_c$, then the number of particles that remain after stabilization is $2n+n+Z+O(1)$.  In this case, we have $\rho_n(\lambda) = \frac{3n+Z}{2n+1} + o(1) \to \frac53+\frac13 e^{-3\lambda}$ in probability. 
\end{proof}

\section{Sandpiles on the Cayley tree}
\label{sec:tree}

\begin{figure}[b]
\centering
\includegraphics[width=\textwidth]{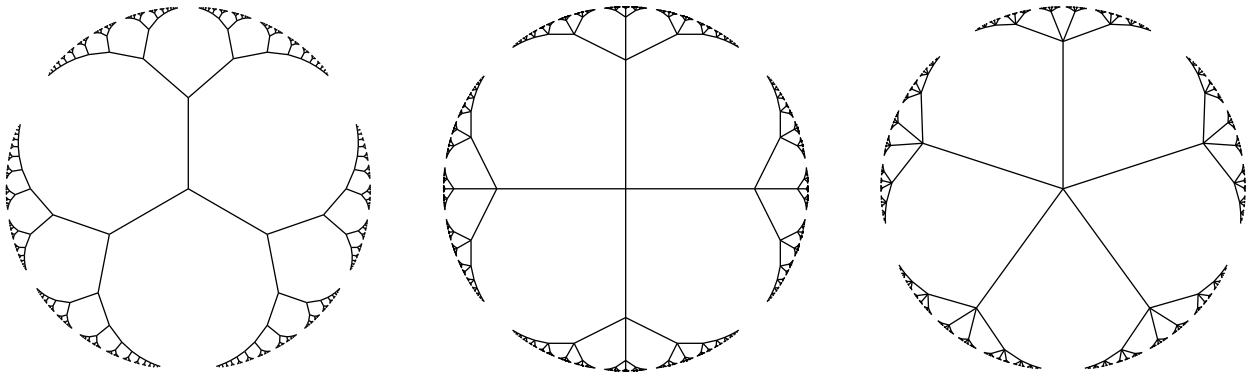}
\caption{The Cayley trees (Bethe lattices) of degree $d=3,4,5$.}
\label{fig:bethe}
\end{figure}

Dhar and Majumdar~\cite{DM} studied the abelian sandpile model on the Cayley tree
(also called the Bethe lattice) with branching factor $q$, which has degree $q+1$.
Implicit in their formulation is that they used wired boundary conditions, i.e., where
all the vertices of the tree at a certain large distance from a
central vertex are glued together and become the sink.  (The other
common boundary condition is free boundary conditions, where all the
vertices at a certain distance from the central vertex become leaves,
and one of them becomes the sink.  The issue of boundary conditions
becomes important for trees, because in any finite subgraph, a constant fraction of vertices are on the boundary.  This is in contrast to $\Z^2$, where free and wired boundary conditions
lead to the same infinite volume limit.  See~\cite{LP:book}.)

The finite regular wired tree $T_{q,n}$ is the ball of radius $n$ in the infinite $(q+1)$-regular tree, with all leaves collapsed to a single vertex $s$.  In $\hat T_{q,n}$ the vertex $s$ serves as the sink.  
Maes, Redig and Saada~\cite{MRS} show that the stationary measure on recurrent sandpiles on $\hat T_{q,n}$ has an infinite-volume limit, which is a measure on sandpiles on the infinite tree.  Denoting this measure by $\Pr_q$, if $h$ denotes the number of particles at a single site far from the boundary, then we have \cite{DM}

	\[ {\textstyle \Pr_q}[h=i] = \frac{1}{(q^2-1)q^q} \sum_{m=0}^i \binom{q+1}{m} (q-1)^{q+1-m} . \]
From this formula we see that the stationary density is
 	\[ \zeta_s = \E_q[h] = \frac{q+1}{2}. \]
For 3-regular, 4-regular, and 5-regular trees, these values are summarized below:

\medskip

\begin{center}
\small
\begin{tabular}{|cc|c|ccccc|}
\hline
\multicolumn{2}{|c|}{tree} & \multirow{2}{*}{$\E_q[h]$} & \multicolumn{5}{c|}{distribution of height $h$ of sand} \\
\cline{1-2}
\cline{4-8}
$q$&degree&& $\Pr_q[h=0]$ & $\Pr_q[h=1]$ & $\Pr_q[h=2]$ & $\Pr_q[h=3]$ & $\Pr_q[h=4]$ \\
\hline
2 & 3 & $3/2$ & 1/12 & 4/12 & 7/12 &&\\
3 & 4 & 2 & 2/27 & 2/9 & 1/3 & 10/27 &\\
4 & 5 & 5/2 & 81/1280 &  27/160 &  153/640 &  21/80  & 341/1280\\
\hline
\end{tabular}
\end{center}

\medskip

Large-scale simulations on $T_{q,n}$ are rather impractical because the vast majority of vertices are near the boundary.  Consequently, each simulation run produces only a small amount of usable data from vertices near the center.    

To experimentally measure $\zeta_c$ for the Cayley trees, we generated large random regular graphs $G_{q,n}$, and used these as finite approximations of the infinite Cayley tree.Ê 
 Ê 
We used the following procedure to generate random connected bipartite multigraphs of degree $q+1$ on $n$ vertices ($n$ even).  
Let $M_0$ be the set of edges $(i,i+1)$ for $i=1,3,5,\ldots,n-1$.  Then take the union of $M_0$ with $q$ additional i.i.d.\ perfect matchings $M_1,\ldots,M_q$ between odd and even vertices.  Each $M_j$ is chosen uniformly among all odd-even perfect matchings whose union with $M_0$ is an $n$-cycle.  

Most vertices of $G_n$
will not be contained in any cycle smaller than $\log_q n + O(1)$ (see e.g., \cite{bollobas}), so these graphs are locally tree-like.  For this reason, we believe that as $n \to \infty$ the threshold density $\Lambda_c(G_n)$ will be concentrated at the threshold density of the infinite tree. 

Since the choice of multigraph affects the estimate of $\zeta_c$, we generated a new independent random multigraph for each trial.  The results for random regular graphs of degree $3$, $4$ and $5$ are summarized in Tables~~\ref{T3}, \ref{T4}, and \ref{T5}.
We find that for the 5-regular tree, the threshold density
is about 2.511 rather than 2.5, for the 4-regular tree the threshold
density is very close to but decidedly larger than 2, while for the
3-regular tree the threshold density is extremely close to 1.5, with
a discrepancy that we were unable to measure.  However, for the 3-regular tree there is
a measurable discrepancy (about $2\times 10^{-6}$) in the probability
that a site has no particles.

\begin{table}[bh]
\begin{center}
\small
\begin{tabular}{|r|r|c|c|c|c|c|}
\hline
\multirow{2}{*}{$n$\ \ } & \multirow{2}{*}{\#samples} & \multirow{2}{*}{$\E[h]$} &
 \multicolumn{3}{c|}{distribution of height $h$ of sand} & \#topplings \\
\cline{4-6}
 &  & & $\Pr[h=0]$ & $\Pr[h=1]$ & $\Pr[h=2]$ &  $\div n\log^{1/2}n$ \\
\hline
1048576 & 2097152 & 1.5004315 & 0.0833326 & 0.332903 & 0.583764 & 1.263145 \\
2097152 & 1048576 & 1.5003054 & 0.0833321 & 0.333031 & 0.583637 & 1.258046 \\
4194304 & 524288 & 1.5002161 & 0.0833314 & 0.333121 & 0.583548 & 1.253092 \\
8388608 & 262144 & 1.5001528 & 0.0833311 & 0.333185 & 0.583484 & 1.247642 \\
16777216 & 131072 & 1.5001081 & 0.0833311 & 0.333230 & 0.583439 & 1.242359 \\
33554432 & 65536 & 1.5000765 & 0.0833307 & 0.333262 & 0.583407 & 1.237317 \\
67108864 & 32768 & 1.5000540 & 0.0833307 & 0.333285 & 0.583385 & 1.232398 \\
134217728 & 16384 & 1.5000382 & 0.0833307 & 0.333300 & 0.583369 & 1.227548 \\
268435456 & 8192 & 1.5000269 & 0.0833308 & 0.333311 & 0.583358 & 1.222371 \\
536870912 & 4096 & 1.5000191 & 0.0833308 & 0.333319 & 0.583350 & 1.214431 \\
1073741824 & 2048 & 1.5000136 & 0.0833307 & 0.333325 & 0.583344 & 1.212751 \\
\hline
\multicolumn{2}{|c|}{$\infty$ (stationary)}
                      & 1.5 & 0.0833333 & 0.333333 & 0.583333 & \multicolumn{1}{c}{} \\
\cline{1-6}
\end{tabular}
\end{center}
\caption{Data for the fixed-energy sandpile on a pseudorandom 3-regular graph on $n$ nodes.  Each estimate of $\EE[h]$ has standard deviation less than $7 \cdot 10^{-8}$, and each estimate of the marginals $\Pr[h=i]$ has standard deviation less than $3 \cdot 10^{-7}$. 
  The data for $\E[h]$ appears to fit $3/2+\text{const}/\sqrt{n}$ very well,
  and extrapolating to $n\to\infty$ it appears that $\E[h]\to 1.500000$ to six decimal places.
  However, apparently $\Pr[h=0] \to 0.083331 < 1/12$.
}
\label{T3}
\end{table}

\begin{table}
\begin{center}
\small
\begin{tabular}{|r|r|c|c|c|c|c|c|}
\hline
\multirow{2}{*}{$n$\ \ } & \multirow{2}{*}{\#samples} & \multirow{2}{*}{$\E[h]$} &
 \multicolumn{4}{c|}{distribution of height $h$ of sand} & \#topplings \\
\cline{4-7}
 &  & & $\Pr[h=0]$ & $\Pr[h=1]$ & $\Pr[h=2]$ & $\Pr[h=3]$ &  $\div n\log^{1/2}n$ \\
\hline
1048576& 2097152& 2.001109& 0.073884& 0.221887& 0.333466& 0.370763& 0.623322\\
2097152& 1048576& 2.000853& 0.073881& 0.221978& 0.333547& 0.370593& 0.618848\\
4194304& 524288& 2.000688& 0.073880& 0.222037& 0.333599& 0.370484& 0.620894\\
8388608& 262144& 2.000584& 0.073878& 0.222075& 0.333631& 0.370416& 0.631324\\
16777216& 131072& 2.000518& 0.073877& 0.222100& 0.333651& 0.370372& 0.649328\\
33554432& 65536& 2.000477& 0.073877& 0.222114& 0.333664& 0.370345& 0.670838\\
67108864& 32768& 2.000451& 0.073877& 0.222123& 0.333673& 0.370328& 0.691040\\
134217728& 16384& 2.000434& 0.073876& 0.222130& 0.333678& 0.370316& 0.699706\\
268435456& 8192& 2.000424& 0.073876& 0.222134& 0.333681& 0.370310& 0.695065\\
536870912& 4096& 2.000417& 0.073876& 0.222136& 0.333683& 0.370305& 0.684507\\
1073741824& 2048& 2.000413& 0.073876& 0.222138& 0.333684& 0.370303& 0.673061\\
\hline
\multicolumn{2}{|c|}{$\infty$ (stationary)}
                      & 2 & 0.074074 & 0.222222 & 0.333333 & 0.370370 & \multicolumn{1}{c}{} \\
\cline{1-7}
\end{tabular}
\end{center}
\caption{Data for the  fixed-energy sandpile on a pseudorandom 4-regular graph on $n$ nodes.  Each estimate of $\E[h]$ and of the marginals $\Pr[h=i]$ has standard deviation less than $3 \cdot 10^{-7}$.
}
\label{T4}
\end{table}

\begin{table}
\begin{center}
\small
\begin{tabular}{|r|r|c|c|c|c|c|c|c|}
\hline
\multirow{2}{*}{$n$\ \ } & \multirow{2}{*}{\#samples} & \multirow{2}{*}{$\E[h]$} &
 \multicolumn{5}{c|}{distribution of height $h$ of sand} & \#topplings \\
\cline{4-8}
 &  & & $\Pr[h=0]$ & $\Pr[h=1]$ & $\Pr[h=2]$ & $\Pr[h=3]$ & $\Pr[h=4]$ &  $\div n$ \\
\hline
1048576 & 1048576 & 2.512106 & 0.062271 & 0.166547 & 0.237230 & 0.264711 & 0.269242 & 1.666086 \\

2097152 & 524288 & 2.511947 & 0.062269 & 0.166579 & 0.237256 & 0.264727 & 0.269169 & 1.666244 \\

4194304 & 262144 & 2.511847 & 0.062268 & 0.166599 & 0.237272 & 0.264737 & 0.269123 & 1.666404 \\

8388608 & 131072 & 2.511781 & 0.062267 & 0.166613 & 0.237283 & 0.264743 & 0.269093 & 1.666589 \\

16777216 & 65536 & 2.511743 & 0.062267 & 0.166621 & 0.237289 & 0.264748 & 0.269075 & 1.667322 \\

33554432 & 65536 & 2.511716 & 0.062267 & 0.166627 & 0.237293 & 0.264750 & 0.269063 & 1.668196 \\

67108864 & 32768 & 2.511700 & 0.062267 & 0.166630 & 0.237296 & 0.264752 & 0.269056 & 1.669392 \\

134217728 & 16384 & 2.511689 & 0.062267 & 0.166632 & 0.237297 & 0.264755 & 0.269050 & 1.671613 \\

268435456 & 8192 & 2.511683 & 0.062266 & 0.166634 & 0.237299 & 0.264753 & 0.269048 & 1.675479 \\

536870912 & 4096 & 2.511680 & 0.062267 & 0.166633 & 0.237300 & 0.264755 & 0.269045 & 1.677092 \\

1073741824 & 2048 & 2.511677 & 0.062266 & 0.166634 & 0.237300 & 0.264755 & 0.269044 & 1.688093 \\
\hline
\multicolumn{2}{|c|}{$\infty$ (stationary)}
                      & 2.5 & 0.063281
                      & 0.168750
                      & 0.239063
                      & 0.262500
                      & 0.266406
                      & \multicolumn{1}{c}{} \\
\cline{1-8}
\end{tabular}
\end{center}
\caption{Data for the fixed-energy sandpile on a pseudorandom 5-regular graph on $n$ nodes.  Each estimate of $\EE[h]$ and of the marginals $\Pr[h=i]$ has standard deviation less than $2 \cdot 10^{-6}$.
}
\label{T5}
\end{table}

\clearpage

\section{Sandpiles on the ladder graph}

\begin{figure}[b]
\centering
\includegraphics[width=0.8\textwidth]{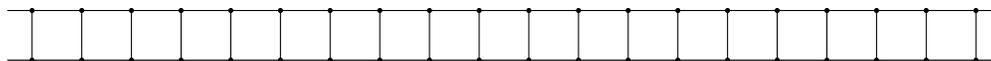}
\caption{The ladder graph.}
\label{fig:ladder}
\end{figure}

The examples in previous sections suggest that the density conjecture can fail for (at least) two distinct reasons: local toppling invariants, and boundary effects.  A \emph{toppling invariant\/} for a graph $G$ is a function $f$ defined on sandpile configurations on $G$
which is unchanged by performing topplings; that is
	\[ f(\sigma) = f(\sigma+\Delta_x) \]
for any sandpile $\sigma$ and any column vector $\Delta_x$ of the Laplacian of $G$.
Examples we have seen are
	\[ f(\sigma) = \sigma(x) \; \mod{2} \]
where $x$ is any vertex of the bracelet graph $B_n$; and
	\[ f(\sigma) = \sigma(x_1) - \sigma(x_2) \; \mod{3} \]
where $x_1,x_2$ are the two vertices comprising any petal on the flower graph $F_n$.  Both of these toppling invariants are \emph{local\/} in the sense that they depend only on a bounded number of vertices as $n \to \infty$.

The Cayley tree has no local toppling invariants, but the large number of sinks, comparable to the total number of vertices, produce a large boundary effect.  The density conjecture fails even more dramatically on the complete graph (\tref{completegraphmain}).  One might guess that this is due to the high degree of interconnectedness, which causes boundary effects from the sink to persist as $n\to \infty$.  A good candidate for a graph $G$ satisfying the density conjecture, then, should have
	\begin{itemize}
	\item no local toppling invariants,
	\item most vertices far from the sink.
	\end{itemize}
The best candidate graphs $G$ should be essentially one-dimensional, so that the sink is well insulated from the bulk of the graph, keeping boundary effects to a minimum.  Indeed, the only graph known to satisfy the density conjecture is the infinite path $\Z$.
	
J\'{a}rai and Lyons \cite{JL} study sandpiles on graphs of the form $G
\times P_n$, where $G$ is a finite connected graph and $P_n$ is the
path of length $n$, with the endpoints serving as sinks.  The simplest such graphs that are not paths are
obtained when $G=P_1$ has two vertices and one edge.  These graphs are
a good candidate for $\zeta_c = \zeta_s$, for the reasons described above.
Nevertheless, we find that while $\zeta_c$ and $\zeta_s$ are
very close, they appear to be different.

First we calculate $\zeta_s$.  Jarai and Lyons \cite[section 5]{JL}
define recurrent configurations as Markov chains on the state space
\[ X = \left\{(3,3), (3,2), (2,3), (3,1), (1,3), \overline{(3,2)},
  \overline{(2,3)} \right\} \]
describing the possible transitions
from one rung of the ladder to the next.  States $(i,j)$ and
$\overline{(i,j)}$ both represent rungs whose left vertex has $i-1$
particles and whose right vertex has $j-1$ particles.  The distinction between
states $(3,2)$ and $\overline{(3,2)}$ lies only in which transitions
are allowed.
The adjacency matrix describing the allowable transitions is given by
	\[ A = \left( \begin{array}{ccccccc} 1&1&1&1&1&0&0 \\ 1&1&1&1&1&0&0 \\ 1&1&1&1&1&0&0 \\ 1&0&0&0&0&1&0 \\ 1&0&0&0&0&0&1 \\ 1&0&0&0&0&1&0 \\ 1&0&0&0&0&0&1 \end{array} \right). \]
Its largest eigenvalue is $2+\sqrt{3}$, and the corresponding left and right eigenvectors are
	\begin{align*}
		u &= (1+\sqrt{3}, 1+\sqrt{3}, 1+\sqrt{3}, 1, 1, 1, 1) \\
     		v &= (3+\sqrt{3}, 1+\sqrt{3}, 1+\sqrt{3}, 1+\sqrt{3}, 1+\sqrt{3}, 1, 1)^{T}
     	\end{align*}
By the Parry formula \cite{parry}, the stationary probabilities are given by $p(i) = u_i v_i/Z$, where $Z$ is a normalizing constant. So
	\begin{align*}
		p(3,3) &= (1+\sqrt{3})(3+\sqrt{3}) / Z \\
     		p(2,3)=p(3,2) &= (1+\sqrt{3})^2 / Z \\
     		p(1,3)=p(3,1) &= (1+\sqrt{3}) / Z \\
     		p(\overline{2,3}) = p(\overline{3,2}) &= 1 / Z
     	\end{align*}
where
	\[ Z= (1+\sqrt{3})(3+\sqrt{3}) + 2(1+\sqrt{3})^2 + 2(1+\sqrt{3}) + 2.
			\]
Thus we find that for the ladder graph in stationarity, the number $h$ of particles at a site satisfies
\begin{align*}
 \Pr[h=0] &=\textstyle -\frac12+\frac{\sqrt{3}}{3} = 0.0773503\dots \\
 \Pr[h=1] &=\textstyle \frac54-\frac{7\sqrt{3}}{12} = 0.2396370\dots \\
 \Pr[h=2] &=\textstyle \frac14+\frac{\sqrt{3}}{4} = 0.6830127\dots\\
 \zeta_s = \E[h] &=\textstyle \frac74 - \frac{\sqrt{3}}{12} = 1.60566243\dots
 \end{align*}
In contrast, the threshold density for ladders appears to be about 1.6082.  Table~\ref{table:ladder-density} summarizes simulation data on finite $2\times n$ ladders.

\begin{table}
\label{table:ladder-density}
\begin{center}
\small
\begin{tabular}{|r|r|c|c|c|c|c|}
\hline
\multirow{2}{*}{$n$\ \ } & \multirow{2}{*}{\#samples} & \multirow{2}{*}{$\E[h]$} &
 \multicolumn{3}{c|}{distribution of height $h$ of sand} & \#topplings \\
\cline{4-6}
 &  & & $\Pr[h=0]$ & $\Pr[h=1]$ & $\Pr[h=2]$ &  $\div n^{5/2}$ \\
\hline
256 & 4194304 & 1.60567 & 0.07695 & 0.24043 & 0.68262 & 0.094773 \\
512 & 2097152 & 1.60693 & 0.07656 & 0.23996 & 0.68349 & 0.095366 \\
1024 & 1048576 & 1.60757 & 0.07636 & 0.23970 & 0.68393 & 0.095864 \\
2048 & 524288 & 1.60788 & 0.07626 & 0.23960 & 0.68414 & 0.096316 \\
4096 & 262144 & 1.60805 & 0.07621 & 0.23952 & 0.68426 & 0.096545 \\
8192 & 131072 & 1.60814 & 0.07618 & 0.23950 & 0.68432 & 0.096753 \\
16384 & 65536 & 1.60816 & 0.07618 & 0.23949 & 0.68434 & 0.097113 \\
32768 & 32768 & 1.60818 & 0.07617 & 0.23948 & 0.68435 & 0.096944 \\
65536 & 16384 & 1.60820 & 0.07616 & 0.23948 & 0.68436 & 0.097342 \\
131072 & 8192 & 1.60820 & 0.07617 & 0.23946 & 0.68437 & 0.097648 \\
262144 & 4096 & 1.60821 & 0.07615 & 0.23949 & 0.68436 & 0.096158 \\
\hline
\multicolumn{2}{|c|}{$\infty$ (stationary)}
     & 1.60566  & 0.07735    & 0.23964   & 0.68301   & \multicolumn{1}{c}{} \\
\cline{1-6}
\end{tabular}
\end{center}
\caption{Data for the fixed-energy sandpile on $2\times n$ ladder graphs.  Each estimate of $\EE[h]$ and of the marginals $\Pr[h=i]$ has a standard deviation smaller than $10^{-5}$.   To four decimal places, the threshold density $\zeta_c$ equals $1.6082$, which exceeds the stationary density $\zeta_s =  7/4-\sqrt{3}/12 = 1.6057$.   The total number of topplings appears to scale as $n^{5/2}$.
}
\end{table}

\newpage

\bibliographystyle{halpha}
\bibliography{soc}

\end{document}